\documentclass[12pt]{article}
\usepackage{amsmath,amsopn,amssymb,amsthm,amsfonts}
\usepackage[T2A]{fontenc}
\usepackage[cp1251]{inputenc}
\usepackage{longtable}

\newtheorem{theorem}{Theorem}[section]

\newtheorem{definition}[theorem]{Definition}

\newcommand{\w}{\omega}
\newcommand\g{{\mathfrak g}}

\newcommand\h{{\mathfrak h}}

\begin{document}

{\bf \large
\centerline{Smolentsev~N.~K., Sokolova~A.~Yu.}

\vspace{3mm}
\centerline{Left-invariant K\"{a}hler and semi-para-K\"{a}hler structures }
\centerline{on some six-dimensional unsolvable Lie groups}
\vspace{3mm}
}

\begin{abstract}
In this article studies questions about the existence of left-invariant K\"{a}hler and semi-para-K\"{a}hler structures on six-dimensional unsolvable Lie groups whose Lie algebras are semidirect products. According to the classification results, there are four such Lie algebras.
It is shown that one of these four Lie algebras admits left-invariant K\"{a}hler metrics, and the other three admit left-invariant semi-para-K\"{a}hler and semi-K\"{a}hler structures.
The paper also shows that a six-dimensional symplectic Lie algebra must be solvable except in one case.
\end{abstract}

\section{Introduction} \label{Intro}
Left-invariant structures on a Lie group are determined by their values on the Lie algebra. Therefore, when studying them, we usually deal with structures on the Lie algebra of a Lie group. In this sense, we will further refer, for example, about symplectic or K\"{a}hler Lie algebra, meaning a left-invariant symplectic or K\"{a}hler structure on the appropriate Lie group.

It is well known that solvable Lie algebras can admit symplectic structures \cite{Campoamor}, \cite{Goze-Khakim}. In the work of Chu B.-Y. \cite{Chu} it was shown that a four-dimensional symplectic Lie algebra must be solvable. This is not the case in the six-dimensional case, and Chu gave an example of a six-dimensional symplectic but unsolvable Lie algebra.
In the same work it was shown that semisimple Lie algebras cannot be symplectic. If a Lie algebra $\g$ is unsolvable, then it has a L\'{e}vy-Maltsev decomposition $\g = N \otimes S$ in the form of a direct sum of the radical $N$ and a solvable subalgebra $S$.
It was shown in \cite{Chu} if a Lie algebra $\g$ has a L\'{e}vy-Malcev decomposition in the form of a direct product $\g = N\times S$, then there are no symplectic structures on it.
Thus, the question remains about the existence of symplectic structures only on unsolvable six-dimensional Lie algebras, for which the L\'{e}vy-Maltsev expansion is a semidirect product $\g = N\ltimes S$.

According to the classification \cite{Basarab}, \cite{Turk} of six-dimensional unsolvable Lie algebras, which are semidirect products, there are four classes of such Lie algebras.
These are the Lie algebras $A_{3.5}\ltimes sl(2,\mathbb{R})$, $A_{3.3}\ltimes sl(2,\mathbb{R})$, $A_{3.1} \ltimes sl(2,\mathbb{R})$ and $A_{3.1} \ltimes so(3)$. Here $A_{3.1}$, $A_{3.3}$, $A_{3.5}$ are three-dimensional solvable Lie algebras with basis $\{e_1,e_2,e_3\}$ and with the following non-zero Lie brackets: $A_{3.1}$ – commutative algebra; $A_{3.3}$: $[e_1,e_2] = e_3$; $A_{3.5}$: $[e_1,e_2] = e_2$, $[e_1,e_3] = e_3$.
Thus, $A_{3.1} = \mathbb{R}^3$, $A_{3.3}$ is three-dimensional Heisenberg algebra $\h_3$ and $A_{3.5}$ is a three-dimensional solvable and non-unimodular Lie algebra.
The Lie algebras $so(3)$ and $sl(2,\mathbb{R})$ are the usual Lie algebras of skew-symmetric matrices of order 3 and matrices of order 2 with zero trace. On these Lie algebras we choose the usual basis $\{e_4,e_5,e_6\}$.
Thus, there are \cite{Basarab} four six-dimensional undecidable Lie algebras, which are semidirect products.
They are defined by the following Lie brackets:


\begin{enumerate}
  \item $A_{3.5}\ltimes sl(2,\mathbb{R})$ : \ $[e_1,e_2] =e_2$, $[e_1,e_3] = e_3$, $[e_4,e_5] = 2e_5$, $[e_4,e_6] = -2e_6$, $[e_5,e_6] = e_4$, $[e_2,e_4] = e_2$, $[e_2,e_5] = e_3$, $[e_3,e_4] =-e_3$, $[e_3,e_6] = e_2$.
  \item $\h_3\ltimes sl(2,\mathbb{R})$ : \ $[e_1,e_2] = e_3$, $[e_4,e_5] = 2e_5$, $[e_4,e_6] = -2e_6$, $[e_5,e_6] = e_4$, $[e_4,e_1] = e_1$, $[e_5,e_2] = e_1$, $[e_6,e_1] = e_2$, $[e_4,e_2] = -e_2$.
  \item $\mathbb{R}^3 \ltimes sl(2,\mathbb{R})$ :\  $[e_4,e_5] = 2e_5$, $[e_4,e_6] = -2e_6$, $[e_5,e_6] = e_4$, $[e_4,e_1] = 2e_1$, $[e_5,e_2] = 2e_1$, $[e_6,e_1] = e_2$, $[e_4,e_3] = -2e_3$, $[e_5,e_3] = e_2$, $[e_6,e_2] = 2e_3$.
  \item $\mathbb{R}^3 \ltimes so(3)$ : \  $[e_4,e_5] = e_6$, $[e_4,e_6] = -e_5$, $[e_5,e_6] = e_4$, $[e_4,e_2] = e_3$, $[e_5,e_1] = -e_3$, $[e_6,e_1] = e_2$, $[e_4,e_3] = -e_2$, $[e_5,e_3] = e_1$, $[e_6,e_2] = -e_1$.
\end{enumerate}

Let $G_1$, $G_2$, $G_3$ and $G_4$ be Lie groups with Lie algebras of the  list above.

The question of left-invariant symplectic structures on these Lie groups (i.e., symplectic structures of Lie algebras) was considered in the work of the authors \cite{Smolen}.
It is shown that only one of the four Lie algebras admits symplectic structures. Thus, Chu's result on six-dimensional unsolvable Lie algebras can be refined as follows: {\it a six-dimensional symplectic Lie algebra must be solvable except in one case $A_{3.5}\ltimes sl(2,\mathbb{R})$}. Therefore, this exceptional unsolvable symplectic Lie algebra $A_{3.5}\ltimes sl(2,\mathbb{R})$ deserves more detailed study.

This work is a extension of the paper \cite{Smolen}.
We will consider questions about the existence of left-invariant K\"{a}hler, semi-K\"{a}hler and semi-para-K\"{a}hler structures on the above four six-dimensional unsolvable Lie groups.
It will be shown that only on one of the four Lie groups, namely on $G_1$, there exist left-invariant K\"{a}hler structures and even with Einstein pseudo-Riemannian metrics.
On the other three Lie groups there exist left-invariant semi-para-K\"{a}hler and semi-K\"{a}hler structures with integrable complex or paracomplex structures.

\section{Preliminaries} \label{Pre}
Let us recall the basic concepts used in this work.
An almost complex structure on a $2n$-dimensional manifold $M$ is a field $J$ of endomorphisms $J: TM\to TM$ such that $J^2 =-Id$.
The almost complex structure $J$ is called integrable if the Nijenhuis tensor, defined by the equality $N_J(X,Y) = [JX,JY ] -[X,Y] -J[JX,Y] -J[X,JY]$, vanishes, for any vector fields $X, Y$ on $M$. In this case, $J$ defines the structure of a complex manifold on $M$.

We will consider left-invariant almost complex structures on the Lie group $G$, which are definite by the left-invariant endomorphism field $J: TG \to TG$ of the tangent bundle $TG$.
Since such a tensor $J$ is determined by a linear operator on the Lie algebra $\g = T_eG$, we will say that $J$ is an invariant almost complex structure on the Lie algebra $\g$.
In this case, the integrability condition for $J$ is formulated at the Lie algebra level: $N_J(X,Y) = 0$, for any $X,Y\in \g$.
In this case, we will say that $J$ is a complex structure on the Lie algebra $\g$.

A left-invariant symplectic structure $\omega$ on a Lie group $G$ is definite by a closed 2-form of maximum rank. The closedness of the form $\omega$ is equivalent to the following condition on the Lie algebra $\g$: $\omega([X,Y],Z) -\omega([X,Z],Y) + \omega([Y,Z],X) = 0$, $\forall X,Y,Z\in \g$. In this case, the Lie algebra $\g$ will be called symplectic.


A left-invariant K\"{a}hler structure on a Lie group $G$ is a triple $(\w, J, g)$ consisting of a left-invariant symplectic form $\w$, a left-invariant complex structure $J$ compatible with the form $\w$: $\w(JX,JY) = \w(X,Y)$, and the (pseudo)Riemannian metric $g$, defined by the formula $g(X,Y) = \w(X,JY)$.

An almost paracomplex structure on a $2n$-dimensional manifold $M$ is the field $P$ of endomorphisms of the tangent bundle $TM$ such that $P^2 = Id$ and the ranks of the eigendistributions $T^\pm M : = {\text ker}(Id \mp P)$ are equal. The almost paracomplex structure $P$ is called to be integrable if the distributions $T^\pm M$ are involutive.
In this case, $P$ is called a paracomplex structure.
The Nijenhuis tensor $N_P$ of an almost paracomplex structure $P$ is defined by the equality $N_P(X,Y) = [X,Y] +[PX,PY] -P[PX,Y] -P[X,PY ]$, for any vector fields $X, Y$ on $M$.
As in the complex case, a para-complex structure $P$ is integrable if and only if $N_P(X,Y) = 0$. A review of the theory of para-complex structures is presented in \cite{Aleks}.

In the case of a left-invariant paracomplex structure $P$ on a Lie group $G$, the eigenspaces subspaces $T_e^\pm G=\g^\pm$ are subalgebras.
Therefore, the paracomplex Lie algebra $\g$ can be represented as a direct sum of two subalgebras:
$$
\g = \g^+ \oplus \g^-.
$$

The almost paracomplex structure $P$ on a Lie group $G$ is said to be
bi-invariant if it is left-invariant and right-invariant simultaneously. This means that $Ad_gJ = J$ for any element $g\in G$. For a connected Lie group, the last condition is equivalent to the following: $\left[ad_X,J\right] =0$, $\forall X\in \g$, or
$$
\left[X, JY\right] = J\left[X,Y\right] .
$$
If an paracomplex structure $P$ on a Lie algebra $\g$ is bi-invariant, then $\g^+ $ and $\g^-$ are ideals of $\g$. The bi-invariant almost paracomplex structure on a Lie group $G$ is always integrable.

A left-invariant para-K\"{a}hler structure on a Lie group $G$ is a triple $(\w,P,g)$ consisting of a left-invariant symplectic form $\w$, a left-invariant paracomplex structure $P$ compatible with the form $\w$: $\w(PX,PY) =-\w(X,Y)$, and the (pseudo)Riemannian metric $g$, defined by the formula $g(X,Y) = \w(X,PY)$.

On the last three Lie algebras of the list presented in the introduction, the condition $d\w = 0$ leads to degeneracy of the form $\w$.
At the same time, $d(\w^3) = 3\w\wedge\w\wedge d\w = 0$ for any 2-form $\w$ on a six-dimensional manifold.
An intermediate property would be $d(\w^2) = 2\w\wedge d\w = 0$.
Therefore, we reduce the closedness property $d\w = 0$ and require that the following property hold:
$$
\w\wedge d\w = 0.
$$
Note that in the case of a left-invariant almost Hermitian structure on a Lie group $G$ of dimension $2n$, the property $d(\w^{n-1}) = 0$ of the fundamental form $\w$ defines a class of semi-K\"{a}hler manifolds according to the Gray–Harvella classification \cite{Gray-Harv}, i.e. such that $\delta \w = 0$.
In our case of pseudo-Hermitian metrics on Lie groups, we will by analogy call such manifolds semi-K\"{a}hler.

\begin{definition}\label{1}
A left-invariant almost Hermitian (almost para-Hermitian) structure $(G,g,\w,J)$ of dimension $2n$, whose fundamental form $\w$ has the property $d(\w^{n-1}) = 0$ is called semi-K\"{a}hler (semi-para-K\"{a}hler).
\end{definition}

The Lie algebras $\h_3\ltimes sl(2,\mathbb{R})$, $\mathbb{R}^3\ltimes sl(2,\mathbb{R})$ and $\mathbb{R}^3\ltimes so(3)$ do not have symplectic forms.
For a non-degenerate 2-form $\w$ its exterior derivative $d\w$ is not equal to zero. In this case, $d\w$ can be non-degenerate as a 3-form. The concept of non-degeneracy (stability) for a 3-form $\Omega$ on a six-dimensional space was defined in Hitchin’s work \cite{Hitch}.
For the 3-form $\Omega$, Hitchin constructed a linear operator $K_\Omega$ whose square is proportional to the identity operator $Id$. Let us recall the basic constructions of Hitchin.

Let $V$ be a 6-dimensional real vector space and $\mu$ a volume form on $V$. Let $\Omega\in \Lambda^3 V^*$ and $X\in V$, then $i_X\Omega\in \Lambda^2 V^*$ and $i_X\Omega\wedge\Omega\in  \Lambda^5 V^*$.
The natural pairing by the exterior product $V^*\otimes \Lambda^5 V^* \to \Lambda^6 V^* \cong \mathbb{R}\mu$ defines an isomorphism $A:\Lambda^5 V^* \cong  V$ and using this, Hitchin defined the linear transformation $K_\Omega : V \to V$ by the following formula:
$$
K_\Omega(X) = A(i_X\Omega\wedge\Omega).
$$
In other words, $i_{K_\Omega(X)}\mu = i_X\Omega\wedge\Omega$.
The operator $K_\Omega$ has the following properties: $\mathrm{trace}(K_\Omega) = 0$ and $K_\Omega^2 = \lambda(\Omega)Id$. Therefore, in the case $\lambda(\Omega)< 0$, we obtain the structure $J_\Omega$ of a complex vector space on the space $V$:
$$
J_\Omega=\frac{1}{\sqrt{-\lambda(\Omega)}} K_\Omega,
$$
and if $\lambda(\Omega)> 0$, then we obtain the para-complex structure $P_\Omega$, i.e., $P_\Omega^2 = Id$, $P_\Omega\ne 1$ by a similar formula:
$$
P_\Omega=\frac{1}{\sqrt{\lambda(\Omega)}} K_\Omega.
$$

Thus, the operator $K_{d\w}$ can define either an almost complex or almost paracomplex structure on six-dimensional Lie algebra $\g$ when $d\w$ is non-degenerate and $\mu = \w^n$.

Let $\nabla$ be the Levi-Civita connection corresponding to the (pseudo)Ri\-eman\-ni\-an metric $g$.
It is determined from a six-term formula \cite{KN}, which for left-invariant vector fields $X,Y,Z$ on a Lie group takes the form: $2g(\nabla_XY,Z) = g([X,Y],Z) + g([Z,X],Y ) + g(X,[Z,Y])$.
The curvature tensor is defined by the formula $R(X,Y) = [\nabla_X, \nabla_Y] - \nabla_{[X,Y]}$.
The Ricci tensor $Ric(X,Y)$ is the convolution of the curvature tensor along the first and fourth (upper) indices.
The Ricci operator $RIC$ is defined by the formula $Ric(X,Y) = g(RIC(X),Y)$.
To calculate the geometric characteristics of left-invariant structures on Lie groups, the Maple system was used.

\section{Left-invariant structures on some unsolvable Lie groups} \label{Groups}

In this section we consider left-invariant geometric structures on each of the four six-dimensional unsolvable Lie groups whose Lie algebras are semidirect products.
Recall that left-invariant geometric structures on a Lie group are determined by their values on the Lie algebra; therefore, when studying them, only the Lie algebra of the Lie group is usually used.
In this sense, we will refer, for example, about a symplectic or K\"{a}hler Lie algebra, meaning a left-invariant symplectic or K\"{a}hler structure on the corresponding Lie group.

\subsection{Left-invariant K\"{a}hler structures on the group $G_1$.} \label{G1}
The group $G_1$ has a Lie algebra $A_{3.5}\ltimes sl(2,\mathbb{R})$ with commutation relations: $[e_1,e_2] =e_2$, $[e_1,e_3] = e_3$, $[e_4,e_5] = 2e_5$, $[e_4,e_6] = -2e_6$, $[e_5,e_6] = e_4$, $[e_2,e_4] = e_2$, $[e_2,e_5] = e_3$, $[e_3,e_4] =-e_3$, $[e_3,e_6] = e_2$.
This Lie algebra is considered by Chu \cite{Chu} as an example showing that a six-dimensional symplectic Lie algebra need not be solvable. It consists of matrices of the form:
$$
\left(
  \begin{array}{ccc}
   a_1 & a_2 &a_3 \\
   0 & a_4 &a_5 \\
   0 & a_6 & -a_4 \\
  \end{array}
\right)
$$

The general form of a symplectic 2-form $\w$ on the Lie algebra $A_{3.5}\ltimes sl(2,\mathbb{R})$ was obtained in \cite{Smolen}:
\begin{equation} \label{1}
\w=\left(
  \begin{array}{cccccc}
   0 & w_{12} & w_{13} & 0 & 0 & 0 \\
   -w_{12} & 0 & 0 & w_{12} & w_{13} & 0 \\
   -w_{13} & 0 & 0 & -w_{13} & 0 & w_{12} \\
   0 & -w_{12} & w_{13} & 0 & w_{45} & w_{46} \\
   0 & -w_{13} & 0 & -w_{45} & 0 & w_{56} \\
   0 & 0 & -w_{12} & -w_{46} & -w_{56} & 0 \\
   \end{array}
\right)
\end{equation}
with the condition of non-degeneracy, $\det(\w) =w_{12}^4 w_{45}^2 -4w_{12}^3 w_{13} w_{45} w_{56} +2w_{12}^2 w_{13}^2 w_{45} w_{46} +4w_{12}^2 w_{13}^2 w_{56}^2-4w_{12} w_{13}^3 w_{46} w_{56}+ w_{13}^4 w_{46}^2 \ne 0$.

It is practically impossible to find a complex structure in a general form compatible with (\ref{1}).
Therefore, we will use special cases where the symplectic 2-form has the simplest form.
From the expression $\det(\w)$ we see that $w_{12}$ and $w_{13}$ cannot be equal to zero at the same time.

Let us first consider the first case, when $w_{12}=1$ and $w_{13}=0$.
Then $\det(\w) =w_{12}^4 w_{45}^2$ does not depend on $w_{46}$ and $w_{56}$. Therefore, we will consider them equal to zero, $w_{46}=0$ and $w_{56}=0$. Let, in addition, $w_{45}=-1$.

In the second case, we will make a similar choice: $w_{12}=0$ and $w_{13}=1$, $w_{56}=0$ and $w_{45}=0$, $w_{46}=1$.

\subsubsection{Case 1.}
Let us consider the symplectic structure represented by
 in Chu's work \cite{Chu}:
\begin{equation} \label{2}
 \w_1 = e^1 \wedge e^2 + e^2 \wedge e^4 + e^3 \wedge e^6 - e^4 \wedge e^5.
\end{equation}

Let us find complex structures $J$ compatible with 2-form (\ref{2}).
To do this, we take a matrix of almost complex structure $J = (J_i^j)$ and require the compatible condition $\w_1(JX,JY) = \w_1(X,Y)$, the condition $J^2 = -Id$ and the integrability condition $N_J(X,Y) = [JX,JY] -[X,Y] -J[X,JY] - J[JX,Y] = 0$. The following system of equations is obtained for finding the compatible complex structure $J$:
\begin{equation} \label{3}
\left\{
\begin{tabular}{l}
$\w_{kj}J^k_i+\w_{ik}J^k_j=0$, \\
$J_k^i\, J^k_j = -\delta_j^i,$  \\
$J_i^l J_j^m C_{lm}^k-J_i^l J_m^k C_{lj}^m-J_j^l J_m^k C_{il}^m -C_{ij}^k=0$,
 \end{tabular}
\right.
\end{equation}
where $\delta_j^i$ is the identity matrix, $C_{ij}^k$ are the structure constants of the Lie algebra, and the indices vary from 1 to $6$.

This system of algebraic equations is solved using Maple symbolic calculations. First, we solve the simplest compatible condition, then the integrability condition, and finally we check the condition $J^2 = -Id$.

For the Lie algebra $A_{3.5}\ltimes sl(2,\mathbb{R})$, we obtain six solutions to the system of equations  (\ref{3}).
For each case, the associated pseudo-Riemannian metric is defined by the formula $g_{1,J}(X,Y) =\w_1(X,JY)$.
As a result, each solution represents a K\"{a}hler structure $(\w_1, J, g_{1,J})$. Moreover, among all the solutions there is one K\"{a}hler structure with an Einstein metric $g_1$ of scalar curvature $S = 12$.
Such a complex structure and the metric tensor are represented by the following matrices:
\begin{equation} \label{4}
J_1=\left(
  \begin{array}{cccccc}
   0 & 1 & 0 & 0 & -\frac 14 & -1 \\
   -1 & 0 & 0 & 1 & 0 & 0 \\
   0 & 0 & 0 & 0 & 1 & 0 \\
   0 & 0 & 0 & 0 & -\frac 14 & -1 \\
   0 & 0 & -1 & 0 & 0 & 0 \\
   0 & 0 & \frac 14 & 1 & 0 & 0
   \end{array}
\right), \qquad
g_1=\left(
  \begin{array}{cccccc}
    -1 & 0 & 0 & 1 & 0 & 0 \\
    0 & -1 & 0 & 0 & 0 & 0 \\
    0 & 0 & \frac 14 & 1 & 0 & 0 \\
    1 & 0 & 1 & -1 & 0 & 0 \\
    0 & 0 & 0 & 0 & -\frac 14 & -1 \\
    0 & 0 & 0 & 0 & -1 &0 \\
   \end{array}
\right).
\end{equation}

\subsubsection{Case 2.}
Let us consider the second special case of the symplectic forms (1):
\begin{equation} \label{5}
 \w_2 = e^1 \wedge e^3 + e^2 \wedge e^5 - e^3 \wedge e^4 + e^4 \wedge e^6.
\end{equation}

Calculations show that there are five families of complex structures $J$ compatible with form (\ref{5}).
Each solution represents a K\"{a}hler structure $(\w_2, J, g_{2,J})$. Moreover, among all the solutions there is one K\"{a}hler structure $(\w_, J_2, g_{2,J})$ with an Einstein metric $g_2$ of nonzero scalar curvature $S = -12a$, depending on one parameter.
This complex structure and metric tensor are represented by the following matrices:
\begin{equation} \label{6}
J_2=\left(
  \begin{array}{cccccc}
    0 & 0 & -a & 0 & \frac 1a & \frac a4 \\
    0 & 0 & 0 & 0 & 0 & -\frac 1a \\
    \frac 1a & 0 & 0 & \frac 1a & 0 & 0  \\
    0 & 0 & 0 & 0 & -\frac 1a & -\frac a4 \\
    0 & -\frac {a^3}{4} & 0 & a & 0 & 0 \\
    0 & a & 0 & 0 & 0 & 0 \\
  \end{array}
\right), \qquad
g_2=\left(
  \begin{array}{cccccc}
    \frac 1a & 0 & 0 & \frac 1a & 0 & 0 \\
    0 & -\frac {a^3}{4} & 0 & a & 0 & 0 \\
    0 & 0 & a & 0 & 0 & 0 \\
    \frac 1a & a & 0 & \frac 1a & 0 & 0 \\
    0 & 0 & 0 & 0 & 0 & \frac 1a \\
    0 & 0 & 0 & 0 & \frac 1a & \frac a4 \\
  \end{array}
\right).
\end{equation}

\subsubsection{Para-Hermitian structure.}\label{sub-3-G1}
The Lie algebra $A_{3.5}\ltimes sl(2,\mathbb{R})$ has an (integrable) paracomplex structure $P_0$ corresponding to the semidirect product of subalgebras $\g^+ =A_{3.5}$ and  $\g^- =sl(2,\mathbb{R})$.
It has a diagonal matrix $P_0 = \text{diag}\{1,1,1,-1,-1,-1\}$.
In addition, on the Lie algebra $A_{3.5}\ltimes sl(2,\mathbb{R})$ there is a 2-form $\Omega_0 = e^1 \wedge e^4 + e^2 \wedge e^5 + e^3 \wedge e^6$, also corresponding to the decomposition into a semidirect product (here $e^1,\dots, e^6$ is the dual basis).
However, the form $\Omega_0$ is not semi-K\"{a}hler.
Therefore, on the Lie algebra $A_{3.5}\ltimes sl(2,\mathbb{R})$ the para-Hermitian structure $(\Omega_0, P_0, G_0)$ is defined, where the pseudo-Riemannian metric $g_0$ is defined by the formula $G_0(X,Y) = \Omega_0(X,P_0Y)$.
This metric has the following Ricci tensor
$$
Ric=-e^1\odot e^1 -\frac 12 e^5\odot e^5 -\frac 12 e^6\odot e^6 -2e^2\odot e^5 +2e^3\odot e^6 -2e^4\odot e^3 -5e^5\odot e^6
$$
and zero scalar curvature.
Recall that  $e^i \odot e^j = \frac 12(e^i \otimes e^j +e^i\otimes e^j)$ is a symmetric product.

Let us formulate the results obtained in the form of a theorem.

\begin{theorem} \label{G1}
The Lie group $G_1$ with the Lie algebra $A_{3.5}\ltimes sl(2,\mathbb{R})$ admits a multiparameter family of left-invariant symplectic structures (\ref{1}). The group $G_1$ admits left-invariant K\"{a}hler structures, including those that have Einstein metrics (\ref{4}) and (\ref{6}). The group $G_1$ also has a natural left-invariant para-Hermitian structure of zero scalar curvature.
\end{theorem}

\subsection{Left-invariant K\"{a}hler structures on the group $G_2$.} \label{G2}
The group $G_2$ has a Lie algebra $\h_3\ltimes sl(2,\mathbb{R})$ with commutation relations: \ $[e_1,e_2] = e_3$, $[e_4,e_5] = 2e_5$, $[e_4,e_6] = -2e_6$, $[e_5,e_6] = e_4$, $[e_4,e_1] = e_1$, $[e_5,e_2] = e_1$, $[e_6,e_1] = e_2$, $[e_4,e_2] = -e_2$.

There are no left-invariant symplectic structures on this Lie algebra; the closedness condition $d\w = 0$ implies the degeneracy of the 2-form $\w$.
It's easy to show. Let $\w = \w_{ij}e^i\wedge e^j$ be an arbitrary 2-form. The components of the $d\w$ have the following expressions through structure constants:
$$
d\w_{ijk}= -C_{ij}^s \w_{sk}+ C_{ik}^s \w_{sj} -C_{jk}^s \w_{si}.
$$
From the condition $d\w_{ijk}=0$ we have, in particular, $d\w_{134}=-C_{13}^s \w_{s4} +C_{14}^s \w_{s3} -C_{34}^s \w_{s1} =0$.
From the commutation relations we obtain $C_{13}^s=0$, $C_{34}^s=0$ and $C_{14}^1=-1$.  Therefore we have: $d\w_{134} = -\w_{13} = 0$.
The following equalities are obtained in exactly the same way: $\w_{23} = 0$, $\w_{34} = 0$, $\w_{35} = 0$, $\w_{36} = 0$, from which it follows that the form $\w$ is degenerate.

\subsubsection{Semi-para-K\"{a}ehler structure corresponding to the semidirect product.} \label{sub-1-G2}
On the Lie algebra $\h_3\ltimes sl(2,\mathbb{R})$ there is a paracomplex structure $P_0$ corresponding to the semidirect product of subalgebras: $\g^+ =\h_3$ and $\g^- =sl(2,\mathbb{R})$.
It has a diagonal matrix $P_0 = \text{diag}\{1,1,1,-1,-1,-1\}$.
In addition, the Lie algebra $\h_3\ltimes sl(2,\mathbb{R})$ has a natural 2-form $\Omega = ae^1\wedge e^4 + be^2\wedge e^5 + ce^3\wedge e^6$, also corresponding to the semidirect product.
It is easy to see that it is compatible with the operator $P_0$ and is semi-K\"{a}hler under the condition $a = -b$.
Let us choose the two simplest semi-K\"{a}hler 2-forms:
\begin{equation} \label{7}
\Omega_{01} = -e^1\wedge e^4 + e^2\wedge e^5 + e^3\wedge e^6 \ \text{ and } \ \Omega_{02} = e^1\wedge e^4 - e^2\wedge e^5 + e^3\wedge e^6.
\end{equation}
Thus, on the Lie algebra $\h_3\ltimes sl(2,\mathbb{R})$ two semi-para-K\"{a}ehler structures $(\Omega_{0i}, P_0, G_{0i})$, $i = 1,2$, are defined, whose pseudo-Riemannian metrics $G_{0i}(X,Y) = \Omega_{0i}(X,P_0Y)$ have zero scalar curvatures and Ricci tensors:
$$
Ric_1 = -4 e^2\odot e^4 - 4e^4\odot e^4 -8e^5\odot e^6,
$$
$$
Ric_2 = 4 e^2\odot e^4 - 4e^4\odot e^4 -8e^5\odot e^6.
$$

Both of these forms $\Omega_{0i}$ have degenerate exterior differentials in the Hitchin sense.

\subsubsection{General semi-para-K\"{a}ehler structures.} \label{sub_2 G2}
There are no closed non-degenerate 2-forms on the Lie algebra $\h_3\ltimes sl(2,\mathbb{R})$.
Therefore, we will consider semi-K\"{a}hler 2-forms, i.e. non-degenerate forms $\w$ and satisfying the condition $\w \wedge d\w = 0$.

Let $\w = \w_{ij}e^i\wedge e^j$ be an arbitrary 2-form.
The condition $\w \wedge d\w = 0$ is satisfied for the following parameter values:

\centerline{
$\w_{12}\w_{35} - \w_{13}\w_{25} + \w_{15}\w_{23} = 0, \w_{12}\w_{34} - \w_{13}\w_{24}+ \w_{14}\w_{23} = 0, -\w_{12}\w_{36} + \w_{13}\w_{26} - \w_{16}\w_{23} = 0,$}

\centerline{
$\w_{13}\w_{46} - \w_{14}\w_{36} + \w_{16}\w_{34} + \w_{23}\w_{56} - \w_{25}\w_{36} + \w_{26}\w_{35} = 0,$}

\centerline{
$-\w_{13}\w_{56} + \w_{15}\w_{36} - \w_{16}\w_{35} - \w_{23}\w_{45} + \w_{24}\w_{35} - \w_{25}\w_{34} = 0,$}

\centerline{
$-\w_{34}\w_{56} + \w_{35}\w_{46} - \w_{36}\w_{45} = 0.$}

Under the condition of non-degeneracy of $\w$, there are 7 solutions of this system.
In all cases, the 3-form $d\w$ is non-degenerate in the sense of Hitchin and to each 3-form $d\w$ there corresponds an operator $P_{d\w}$ of a (non-integrable) almost paracomplex structure.
Note that the expressions of the 2-forms $\w$ and the corresponding operators $P_{d\w}$ are quite cumbersome.
Therefore, we consider only those semi-K\"{a}hler 2-forms $\w$ with non-degenerate $d\w$ that are compatible with the paracomplex structure $P_0$, i.e. such that $\w(P_0X, P_0Y) = -\w(X,Y)$.
Then we obtain three semi-K\"{a}hler 2-forms:

\begin{multline} \label{8}
\omega_1 =
e^1\wedge \left(\frac{\w_{16} \w_{34} -\w_{25} \w_{36}+ \w_{26}\w_{35} }{\w_{36}}  e^4 + \frac{\w_{16} \w_{35} -\w_{24} \w_{35}+ \w_{25}\w_{34} }{\w_{36}}  e^5+ \w_{16} e^6 \right) \\ +e^2\wedge (\w_{24} e^4 +\w_{25} e^5 +\w_{26} e^6) +e^3\wedge (\w_{34} e^4 +\w_{35} e^5 +\w_{36} e^6 ),
\end{multline}

\begin{multline} \label{9}
\omega_2 =
e^1\wedge \left(\w_{14} e^4+\w_{15} e^5 +\frac{\w_{24}\w_{35} -\w_{25} \w_{34}}{\w_{35}} e^6 \right)+ \\
e^2\wedge \left(\w_{24} e^4+\w_{25} e^5 -\frac{(\w_{24}\w_{35} -\w_{25} \w_{34})\w_{34}}{\w_{35}^2} e^6\right) + e^3\wedge (\w_{34} e^4 +\w_{35} e^5),
\end{multline}

\begin{equation} \label{10}
\omega_3 =
e^1\wedge (\w_{14} e^4+\w_{15} e^5) +e^2\wedge (\w_{24} e^4 +\w_{26} e^6) +e^3\wedge \w_{34} e^4.
\end{equation}

In each case we obtain a semi-para-K\"{a}hler structure $(\w, P_0, g_\w)$ with the associated metric $g_\w(X,Y) = \w(X,P_0Y)$.
Using the Maple system, the geometric characteristics of the $g_\w$ metric can be easily calculated. In particular, for the last form $\w_3$ the following expressions for the metric tensor and scalar curvature hold:
$$
g_3=
\left(
  \begin{array}{cccccc}
    0&0&0&-\w_{14}&-\w_{15}&0 \\
    0&0&0&-\w_{24}&0&-\w_{26} \\
    0&0&0&-\w_{34}&0&0 \\
    -\w_{14}&-\w_{24}&-\w_{34} &0&0&0 \\
    -\w_{15} &0&0&0&0&0 \\
    0&-\w_{26} &0&0&0&0 \\
   \end{array}
\right), \quad
S_3 =\frac{\w_{34}}{\w_{15} \w_{26} }.
$$

As already noted, each of the above 2-forms has the non-degenerate exterior derivative $d\w$ in the Hitchin sense and, therefore, has the  operator $P_{d\w}$ of almost paracomplex structure.
Let us give the expression for the operator $P_{3}$ corresponding to the 3-form $d\w_3$:
$$
P_{3}(e_1) =e_1 -2\frac{\w_{14}}{\w_{34}}e_3 ,\quad
P_{3}(e_2) =e_2 -2\frac{\w_{24}}{\w_{34}}e_3 ,\quad
P_{3}(e_3) =-e_3,
$$
$$
P_{3}(e_4)= \frac{2\w_{14} \w_{24} +18\w_{15} \w_{26}}{\w_{34}^2}e_3 +e_4,\quad
P_{3}(e_5) =2\frac{\w_{14}}{\w_{34}}e_1 + 6\frac{\w_{15}}{\w_{34}}e_2- \frac{2\w_{14}^2 +6\w_{15} \w_{24}}{\w_{34}^2}e_3 -e_5,
$$
$$
P_{5}(e_6)= 6\frac{\w_{26}}{\w_{34}}e_1 -2\frac{\w_{24}}{\w_{34}}e_2 +\frac{2\w_{24}^2 -6\w_{14} \w_{26}}{\w_{34}^2} e_3 -e_6.
$$

Let us formulate the results obtained in the form of a theorem.

\begin{theorem} \label{Th-G2}
The group $G_2$ with Lie algebra $\h_3\ltimes sl(2,\mathbb{R})$ does not admit left-invariant symplectic structures.
The group $G_2$ has natural left-invariant semi-para-K\"{a}hler structures $(\Omega_{0i},P_0, G_{0i})$, $i = 1,2$, of zero scalar curvature, where the semi-K\"{a}hler 2-forms $\Omega_{0i}$ are represented by formulas (\ref{7}), $P_0 = diag{1,1,1 ,-1,-1,-1}$, and $G_{0i}(X,Y) =\Omega_{0i}(X,P_0Y)$.
The group $G_2$ also admits multiparameter families (\ref{8}) – (\ref{10}) of left-invariant semi-K\"{a}hler 2-forms $\w$
compatible with the paracomplex structure operator $P_0$ and, therefore, it admits multiparameter families of semi-paraK\"{a}hler structures $(\w, P_0, g_\w)$  with associated metrics $g_\w(X,Y) = \w(X,P_0Y)$.
\end{theorem}

\subsection{Left-invariant semi-para-K\"{a}hler structures on the group $G_3$.} \label{G3}
The group $G_3$ has a Lie algebra $\mathbb{R}^3 \ltimes sl(2,\mathbb{R})$ with commutation relations: \ $[e_4,e_5] = 2e_5$, $[e_4,e_6] = -2e_6$, $[e_5,e_6] = e_4$, $[e_4,e_1] = 2e_1$, $[e_5,e_2] = 2e_1$, $[e_6,e_1] = e_2$, $[e_4,e_3] = -2e_3$, $[e_5,e_3] = e_2$, $[e_6,e_2] = 2e_3$.
There are no left-invariant symplectic structures on this Lie algebra; the closedness condition $d\w = 0$ implies the degeneracy of the 2-form $\w$.
This is shown just as simply as in the previous section.

\subsubsection{Semi-para-K\"{a}ehler structure corresponding to a semidirect product.}\label{sub-1-G3}
On the considering Lie algebra there is a paracomplex structure $P_0$ corresponding to the semidirect product of the subalgebras $\g^+=\mathbb{R}^3$ and $\g^-=sl(2,\mathbb{R})$.
It has a diagonal matrix $P_0 = \text{diag}\{1,1,1,-1,-1,-1\}$.
In addition, the Lie algebra $\mathbb{R}^3 \ltimes sl(2,\mathbb{R})$ has a natural 2-form $\Omega = ae^1\wedge e^4 + be^2\wedge e^5 + ce^3\wedge e^6$, corresponding to the semidirect product.
It is easy to see that it is compatible with the operator $P_0$ and is semi-K\"{a}hler under the condition $a = -b$.
Let us choose the two simplest semi-K\"{a}hler 2-forms:
$$
\Omega_{01} = -e^1\wedge e^4 + e^2\wedge e^5 + e^3\wedge e^6 \ \text{ and } \ \Omega_{02} = e^1\wedge e^4 - e^2\wedge e^5 + e^3\wedge e^6,
$$
Then two semi-para-K\"{a}ehler structures $(\Omega_{0i}, P_0, G_{0i})$, $i = 1,2$, are defined on the Lie algebra $\mathbb{R}^3 \ltimes sl(2,\mathbb{R})$.
Their pseudo-Riemannian metrics $G_{0i}(X,Y) =\Omega_{0i}(X,P_0Y)$ have zero scalar curvature and the following Ricci tensors:
$$
Ric_i = -16e^4\odot e^4 -16e^5\odot e^6.
$$
Both 2-forms $\Omega_{0i}$ have degenerate exterior derivative $d\Omega_{0i}$ in the Hitchin sense.

\subsubsection{General semi-para-K\"{a}ehler structures. } \label{sub-2-G3}
There are no closed non-degenerate 2-forms on this Lie algebra.
Let $\w = \w_{ij}e^i\wedge e^j$ be an arbitrary 2-form.
The condition $\w \wedge d\w = 0$ is satisfied for the following parameter values:

\centerline{$w_{12}w_{35} -w_{13}w_{25} +w_{15}w_{23} = 0,\quad  w_{12}w_{34} -w_{13}w_{24} + w_{14}w_{23} = 0,$}
\centerline{$-w_{12}w_{56} -w_{13}w_{45} +w_{14}w_{35} +w_{15}w_{26} -w_{15}w_{34} -w_{16}w_{25} = 0,$}
\centerline{$w_{13}w_{46} -w_{14}w_{36} +w_{16}w_{34} +w_{23}w_{56} -w_{25}w_{36} +w_{26}w_{35} = 0,$}
\centerline{$w_{12}w_{46} -w_{14}w_{26} + w_{16}w_{24} -w_{23}w_{45} + w_{24}w_{35} -w_{25}w_{34} = 0,$ }
\centerline{$-w_{12}w_{36} +w_{13}w_{26} -w_{16}w_{23} = 0.$}

There are 2 solutions to this system of equations with non-degenerate 2-forms $\w$:
\begin{multline*}
    \w_1 = e^1 \wedge (-w_{25}e^4 + w_{15}e^5 -w_{35}e^6) + e^2 \wedge (w_{24}e^4 + w_{25}e^5 + w_{26}e^6) + \\ e^3 \wedge (w_{26}e^4 + w_{35}e^5 -w_{36}e^6) + e^4 \wedge (w_{45}e^5 + w_{46}e^6) +w_{56} e^5 \wedge e^6,
\end{multline*}
\begin{multline*}
    \w_2 = e^1 \wedge (w_{14}e^4 + w_{15}e^5 -w_{35}e^6) + e^2 \wedge (w_{24}e^4 -w_{14}e^5 + w_{34}e^6) + \\ e^3 \wedge (w_{34}e^4 + w_{35}e^5) +
    e^4 \wedge (w_{45}e^5 + w_{46}e^6) +w_{56} e^5 \wedge e^6.
\end{multline*}

Both of these 2-forms have a degenerate exterior derivative $d\w_i$, $i=1,2$.

The Lie algebra $\mathbb{R}^3 \ltimes sl(2,\mathbb{R})$ has the integrable paracomplex structure \\ $P_0 = \text{diag}\{1,1,1,-1,-1,-1\}$, corresponding to the decomposition of the Lie algebra into a semidirect product.
Let us require that the above semi-K\"{a}hler 2-forms $\w_i$ be compatible with the operator $P_0$: $\w_i(P_0X,P_0Y) = -\w_i(X,Y)$.
Then,
\begin{equation} \label{13}
\w_1 =e^1 \wedge (-\w_{25}e^4 + \w_{15}e^5 -\w_{35}e^6) + e^2 \wedge (\w_{24}e^4 + \w_{25}e^5 + \w_{26}e^6) + e^3 \wedge (\w_{26}e^4 + \w_{35}e^5 -\w_{36}e^6),
\end{equation}
\begin{equation} \label{14}
\w_2 =e^1 \wedge (\w_{14}e^4 + \w_{15}e^5 -\w_{35}e^6) + e2 \wedge (\w_{24}e^4 -\w_{14}e^5 + \w_{34}e^6) + e^3 \wedge (\w_{34}e^4 + \w_{35}e^5).
\end{equation}

Let us define the associated metric by the formula $g_i(X,Y) =\w_i(X,P_0 Y)$, $i=1,2$.
Calculations show that every semi-para-K\"{a}hler structure $(\w_i,P_0,g_i)$ has zero scalar curvature.

Let us formulate the results obtained in the form of a theorem.

\begin{theorem} \label{Th-G3}
The group $G_3$ with Lie algebra $\mathbb{R}^3\ltimes sl(2,\mathbb{R})$ does not admit left-invariant symplectic structures.
The group $G_3$ has natural left-invariant semi-para-K\"{a}hler structures
$(\Omega_{0i},P_0, G_{0i})$, $i = 1,2$,
of zero scalar curvature, where the semi-K\"{a}hler 2-forms $\Omega_{0i}$ are represented by formulas (\ref{7}), $P_0 = \text{diag}\{1,1,1,-1,-1,-1\}$, and $G_{0i}(X,Y) =\Omega_{0i}(X,P_0Y)$.
The group $G_3$ also admits multiparameter families (\ref{13}) – (\ref{14}) of left-invariant semi-K\"{a}hler 2-forms $\w_i$ compatible with the paracomplex structure operator $P_0$ and, therefore, it admits multiparameter families of left-invariant semi-paraK\"{a}hler structures $(\w_i, P_0, g_{i})$ with associated metrics $g_{i}(X,Y) =\omega_{i}(X,P_0Y)$ zero scalar curvature.
\end{theorem}

\subsection{Left-invariant semi-K\"{a}hler and semi-para-K\"{a}hler structures on the group $G_4$.} \label{G4}
The group $G_4$ has a Lie algebra $\mathbb{R}^3 \ltimes so(3)$ with commutation relations: \ $[e_4,e_5] = e_6$, $[e_4,e_6] = -e_5$, $[e_5,e_6] = e_4$, $[e_4,e_2] = e_3$, $[e_5,e_1] = -e_3$, $[e_6,e_1] = e_2$, $[e_4,e_3] = -e_2$, $[e_5,e_3] = e_1$, $[e_6,e_2] = -e_1$.

This Lie algebra has a paracomplex structure $P_0$ corresponding to the semidirect product of the subalgebras $\mathbb{R}^3$ and $so(3)$.
It has a diagonal matrix $P_0 = \text{diag}\{1,1,1,-1,-1,-1\}$.
In addition, on the Lie algebra $\mathbb{R}^3 \ltimes so(3)$ there is a natural 2-form $\Omega = ae^1\wedge e^4 + be^2\wedge e^5 + ce^3\wedge e^6$, also corresponding to the semidirect product.
It is easy to see that it is compatible with the operator $P_0$ and is semi-K\"{a}hler for any non-zero values of the parameters $a, b, c$.
Therefore, on the Lie algebra $\mathbb{R}^3 \ltimes so(3)$ the semi-para-K\"{a}ehler structure $(\Omega, P_0, g)$ is defined whose pseudo-Riemannian metric $G(X,Y) = \Omega(X,P_0Y)$ has the Ricci tensor of the form
$$
Ric=\frac{a^2-(b-c)^2}{bc} e^4\odot e^4 -\frac{b^2+(a-c)^2}{ab} e^5\odot e^5 +\frac{c^2-(a-b)^2}{ac} e^6\odot e^6
$$
and zero scalar curvature.

\subsubsection{General semi-para-K\"{a}ehler structures. } \label{G4-1}
There are no closed non-degenerate 2-forms on this Lie algebra.
Let $\w = \w_{ij}e^i\wedge e^j$ be an arbitrary 2-form.
The condition $\w \wedge d\w = 0$ is satisfied for the following parameter values:

\centerline{$-w_{13}w_{56} + w_{15}w_{36} -w_{16}w_{35} + w_{23}w_{46} -w_{24}w_{36} + w_{26}w_{34} = 0,$}
\centerline{$w_{12}w_{36} -w_{13}w_{26} + w_{16}w_{23} = 0,\quad  w_{12}w_{34} -w_{13}w_{24} + w_{14}w_{23} = 0,$}
\centerline{$w_{12}w_{56} -w_{15}w_{26} + w_{16}w_{25} -w_{23}w_{45} + w_{24}w_{35} -w_{25}w_{34} = 0, \quad -w_{12}w_{35} + w_{13}w_{25} -w_{15}w_{23} = 0.$}

There are 2 solutions to this system of equations with non-degenerate 2-forms $\w$:
\begin{multline*}
    \w_1 = e^1 \wedge (w_{14}e^4 + w_{24}e^5 + w_{34}e^6) + e^2 \wedge (w_{24}e^4 + w_{25}e^5 + w_{26}e^6) + \\ e^3 \wedge (w_{34}e^4 + w_{26}e^5 + w_{36}e^6) +
    e^4 \wedge (w_{45}e^5 + w_{46}e^6) +w_{56} e^5 \wedge e^6,
\end{multline*}
\begin{multline*}
    \w_2 = e^1 \wedge (w_{14}e^4 + w_{15}e^5 + w_{34}e^6) + e^2 \wedge (w_{15}e^4 + w_{25}e^5 + w_{35}e^6) + \\ e^3 \wedge (w_{34}e^4 + w_{35}e^5) + e^4 \wedge (w_{45}e^5 + w_{46}e^6) +w_{56} e^5 \wedge e^6.
    \end{multline*}

Let us require that these  semi-K\"{a}hler 2-forms $\w_i$ be compatible with the operator $P_0$: $\w_i(P_0X,P_0Y) = -\w_i(X,Y)$.
Then they take the form:
\begin{equation} \label{15}
\w_1 =e^1 \wedge (\w_{14}e^4 + \w_{24}e^5 + \w_{34}e^6) + e^2 \wedge (\w_{24}e^4 + \w_{25}e^5 + \w_{26}e^6) + e^3 \wedge (\w_{34}e^4 + \w_{26}e^5 + \w_{36}e^6),
\end{equation}
\begin{equation} \label{16}
\w_2 =e^1 \wedge (w_{14}e^4 + w_{15}e^5 + w_{34}e^6) + e^2 \wedge (w_{15}e^4 + w_{25}e^5 + w_{35}e^6) + e^3 \wedge (w_{34}e^4 + w_{35}e^5).
\end{equation}

Let us define the associated metric by the formula $g_i(X,Y) =\w_i(X,P_0 Y)$, $i=1,2$.
Calculations show that every semi-para-K\"{a}hler structure $(\w_i,P_0,g_i)$ has zero scalar curvature.

\subsubsection{Semi-K\"{a}ehler structures.} \label{G4-sub-2}
The Lie algebra $\mathbb{R}^3\ltimes so(3)$ admits integrable complex structures compatible with the para-K\"{a}hler 2-form:
\begin{equation} \label{17}
\Omega_{01} = -e^1\wedge e^4 + e^2\wedge e^5 + e^3\wedge e^6.
\end{equation}
Let us find the matrix $J = (\psi_{ij})$ of this complex structure.

The compatible condition $\Omega_{01}(JX,JY) = \Omega_{01}(X,Y)$ is satisfied for the following values of parameters $\psi_{ij}$:
$\psi_{51} = -\psi_{42}$, $\psi_{61} = -\psi_{43}$, $\psi_{11} = -\psi_{44}$, $\psi_{21} = \psi_{45}$, $\psi_{31} = \psi_{46}$, $\psi_{62} = \psi_{53}$, $\psi_{12} = \psi_{54}$, $\psi_{22} = -\psi_{55}$, $\psi_{32} = -\psi_{56}$, $\psi_{13} = \psi_{64}$, $\psi_{23} = -\psi_{65}$, $\psi_{33} = -\psi_{66}$, $\psi_{24} = -\psi_{15}$, $\psi_{34} = -\psi_{16}$, $\psi_{35} = \psi_{26}$.

Then from the other two equations of system (\ref{3}) we obtain 4, up to the sign $(\pm J)$ of the solution:

\begin{enumerate}
  \item $\psi_{15}=0$, $\psi_{16}=(\psi_{45}^2+\psi_{46}^2)/(\psi_{41} \psi_{45})$, $\psi_{25}=\psi_{45}^2/\psi_{41}$, $\psi_{26} =(\psi_{45} \psi_{46})/\psi_{41}$, $\psi_{36}=\psi_{46}^2/\psi_{41}$, $\psi_{42}=0$, $\psi_{43}=0$, $\psi_{44}=-\psi_{46}/\psi_{45}$, $\psi_{52}=0$, $\psi_{53}=0$, $\psi_{54}=0$, $\psi_{55}=0$, $\psi_{56}=-1$, $\psi_{64}=0$, $\psi_{65}=1$, $\psi_{66}=0$, $\psi_{14}=-(\psi_{45}^2+\psi_{46}^2)/(\psi_{41} \psi_{45} ).$

  \item $\psi_{15}=((\psi_{45} \psi_{44}+\psi_{46} ) \psi_{14})/(\psi_{44}^2+1)$, $\psi_{16}=(\psi_{15} (\psi_{46} \psi_{44}-\psi_{45} ))/(\psi_{45} \psi_{44}+\psi_{46} )$, $\psi_{25}=-(\psi_{15} \psi_{45}^2)/(\psi_{45} \psi_{44}+\psi_{46} )$, $\psi_{26}=-(\psi_{15} \psi_{45} \psi_{46})/(\psi_{45} \psi_{44}+\psi_{46} )$, $\psi_{36}=-(\psi_{15} \psi_{45}^2)/(\psi_{45} \psi_{44}+\psi_{46} )$, $\psi_{41}=-(\psi_{45} \psi_{44}+\psi_{46})/\psi_{15}$, $\psi_{42}=0$, $\psi_{43}=0$, $\psi_{52}=0$, $\psi_{53}=0$, $\psi_{54}=0$, $\psi_{55}=0$, $\psi_{56}=-1$, $\psi_{63}=0$, $\psi_{64}=0$, $\psi_{65}=1$, $\psi_{66}=0.$

  \item $\psi_{14}=-(\psi_{16}^2 \psi_{41}^2+\psi_{46}^2)/(\psi_{41}\psi_{46} )$, $\psi_{15}=\psi_{46}/\psi_{41} $, $\psi_{25}=0$, $\psi_{26}=0$, $\psi_{36}=(\psi_{46}^2)/\psi_{41} $, $\psi_{42}=0$, $\psi_{43}=0$, $\psi_{44}=-(\psi_{16} \psi_{41})/\psi_{46} $, $\psi_{45}=0$, $\psi_{52}=0$, $\psi_{53}=0$, $\psi_{54}=0$, $\psi_{55}=0$, $\psi_{56}=1$, $\psi_{63}=0$, $\psi_{64}=0$, $\psi_{65}=-1$, $\psi_{66}=0.$

  \item $\psi_{14}=-(\psi_{44}^2+1)/\psi_{41} $, $\psi_{15}=0$, $\psi_{16}=0$, $\psi_{25}=0$, $\psi_{26}=0$, $\psi_{36}=0$, $\psi_{42}=0$, $\psi_{43}=0$, $\psi_{45}=0$, $\psi_{46}=0$, $\psi_{52}=0$, $\psi_{53}=0$, $\psi_{54}=0$, $\psi_{55}=0$, $\psi_{56}=-1$, $\psi_{63}=0$, $\psi_{64}=0$, $\psi_{65}=1$, $\psi_{66}=0.$

\end{enumerate}

Each of the complex structures $J$ of this list defines a semi-K\"{a}hler structure $(\Omega_{01}, J, g_J)$, where $g_J(X,Y) = \Omega_{01}(X,JY)$.

We present explicit expressions for the fourth complex structure, the associated metric, the Ricci tensor and scalar curvature:
$$
J_4=\left(
  \begin{array}{cccccc}
    -\psi_{44}&0&0&-\frac{\psi_{44}^2+1}{\psi_{41}} &0&0 \\
    0&0&-1&0&0&0 \\
    0&1&0&0&0&0 \\
    \psi_{41}&0&0&\psi_{44}^2&0&0 \\
    0&0&0&0&0&-1 \\
    0&0&0&0&1&0 \\
  \end{array}
\right),\quad
g_{J_4}=\left(
  \begin{array}{cccccc}
    -\psi_{41}&0&0&-\psi_{44}&0&0 \\
    0&0&0&0&0&-1 \\
    0&0&0&0&1&0 \\
    -\psi_{44}&0&0&\frac{\psi_{44}^2+1}{\psi_{41}}&0&0 \\
    0&0&1&0&0&0 \\
    0&-1&0&0&0&0 \\
  \end{array}
\right),
$$
$$
Ric_4=\left(
  \begin{array}{cccccc}
    -\psi_{41}^2&0&0&-\psi_{41} \psi_{44}&0&0 \\
    0&0&0&0&0&\psi_{41}/2 \\
    0&0&0&0&-\psi_{41}/2&0 \\
    -\psi_{41} \psi_{44}&0&0&-\psi_{44}^2&0&0 \\
    0&0&-\psi_{41}/2&0&2&0 \\
    0&\psi_{41}/2&0&0&0&2 \\
  \end{array}
\right), \quad S_4 = \psi_{41}.
$$

Let us formulate the results obtained in the following form.

\begin{theorem} \label{Th-G4}
The group $G_4$ with the Lie algebra $\mathbb{R}^3\ltimes so(3)$ does not admit left-invariant symplectic structures.
The group $G_4$ has natural left-invariant semi-para-K\"{a}hler structures $(\Omega,P_0, G)$, with zero scalar curvature, where $\Omega = ae^1\wedge e^4 + be^2\wedge e^5 + ce^3\wedge e^6$ and $P_0 = {\text diag}\{1,1,1,-1,-1,-1\}$.
The Lie group $G_4$ also admits multiparameter families (\ref{15}), (\ref{16}) of left-invariant semi-K\"{a}hler 2-forms $\w_i$, $i=1,2$, compatible with the paracomplex structure $P_0$ and, therefore, it admits multiparameter families of left-invariant semi-paraK\"{a}hler structures $(\w_i, P_0, g_i)$ with associated metrics $g_i(X,Y)= \w_i(X,P_0Y)$ of zero scalar curvature.
The group $G_4$ admits four multiparameter families of left-invariant complex structures $J$ compatible with the semi-K\"{a}hler 2-form (\ref{17}) and, therefore, it admits multiparameter families of left-invariant semi-K\"{a}hler structures $(\Omega_{01}, J, g_J)$ with integrable complex structures J, where $g_J(X,Y) = \Omega_{01}(X,JY)$.
\end{theorem}


\begin{thebibliography}{999}

\bibitem{Aleks}
Alekseevsky D.V., Medori C., Tomassini A. Homogeneous para-K\"{a}hler Einstein manifolds. Russ. Math. Surv. 2009, Vol. 64, no. 1, pp. 1–43.

\bibitem{Basarab}	
Basarab-Horwath P., Lahno V., Zhdanov R. The structure of Lie algebras and the classification problem for partial differential equations. Acta Appl. Math. 2001. Vol. 69, pp. 43–94.

\bibitem{Campoamor}	
Campoamor-Stursberg R. Symplectic forms on six-dimensional real solvable Lie algebras I. Algebra Colloquium. 2009. Vol. 16, no. 2, pp. 253-266.

\bibitem{Chu}	
Chu Bon-Yao.  Symplectic homogeneous spaces. Trans. of the Amer. Math. Soc. 1974. Vol. 197, pp. 154–159.

\bibitem{Goze-Khakim}
Goze M., Khakimdjanov Y., Medina A. Symplectic or contact structures on Lie groups. Differential Geom. Appl. 2004, Vol. 21, no. 1, 41--54.

\bibitem{Gray-Harv}
Gray A., Harvella L.M. The sixteen classes of almost Hermitian manifolds and their linear Invariants. Ann. Math. Pura Appl. 1980. Vol. 123, pp. 35–58.

\bibitem{KN}
Kobayashi S. and  Nomizu K. Foundations of Differential Geometry, Vol. 1 and 2. Interscience Publ. New York, London. 1963.

\bibitem{Turk}	
Turkowski P. Low–dimensional real Lie algebras. J. Math. Phys. 1988. Vol. 29, pp. 2139–2144.

\bibitem{Smolen}	
Smolentsev N.K., Sokolova A.Yu.  Parak\"{a}hler and para-Hermitian structures on six-dimensional unsolvable Lie algebras. Izvestiya AltGU. Matematika i mekhanika. 2023. Vol. 4(132), pp. 94–98.

\bibitem{Hitch}
Hitchin N.J. The geometry of three-forms in six dimensions J. Diff. Geom. 2000. Vol. 55, pp. 547–576.

\end{thebibliography}
\end{document}